\documentclass[12pt]{article}

\usepackage[top=30truemm,bottom=30truemm,left=25truemm,right=25truemm]{geometry}

\usepackage{amsmath}
\usepackage{amssymb}
\usepackage{amsthm}
\usepackage{float}
\usepackage{graphicx}
\usepackage{url}

\theoremstyle{plain} 
\newtheorem{thm}{Theorem}[section] 

\theoremstyle{definition} 
\newtheorem{defi}[thm]{Definition}

\newtheorem{prop}[thm]{Proposition}
\newtheorem{rem}[thm]{Remark}

\newcommand{\bs}{\mathbf{s}}

\newcommand{\cF}{\mathcal{F}}

\newcommand{\cT}{\mathcal{T}}

\newcommand{\PSL}{\mathrm{PSL}(2,\mathbb{R})}
\newcommand{\SSS}{\bold{S_\infty^1}}

\newcommand{\HHH}{\bold{H}}

\newcommand{\Teich}{Teichm\"uller\ }

\newcommand{\trop}{\text{trop}}

\newcommand{\tropmfd}{\mathcal{X}_\bs (\mathbb{R}^\trop)}

\newcommand{\trf}{|\text{tr}(f_2)|}
\newcommand{\Lf}{L(f_2)}

\newcommand{\Intr}{\mathrm{Int}}
\newcommand{\ceq}{:=}
\newcommand{\annulus}{A_{1,1}}
\newcommand{\wT}{\widehat{\cT}(\Sigma)}
\newcommand{\wTone}{\widehat{\cT}(\annulus)}
\newcommand{\Fi}{\cF_\infty(\Sigma)}
\newcommand{\wiT}{\widehat{\cT}_\infty(\Sigma)}

\title{The Fenchel--Nielsen twist in terms of the cross ratio coordinates}
\author{Takeru Asaka}
\date{}

\begin{document}

\maketitle

\addtocounter{footnote}{0}

\addtocounter{footnote}{0}

\setlength{\abovedisplayskip}{2pt}
\setlength{\belowdisplayskip}{2pt}

\begin{abstract}
    We calculate the Fenchel--Nielsen twist in the enhanced \Teich space of a marked surface by the cross ratio coordinates.
\end{abstract}


\section{Introduction}
\hspace{20pt} The Fenchel--Nielsen twist is a generalization of the Dehn twist. It is a deformation in the \Teich space defined by cutting a hyperbolic surface along a simple closed geodesic, rotating one side of the cut relative to the other, and attaching these sides. If the scale of the rotation is equal to the length of the simple closed geodesic, then the Fenchel--Nielsen twist is equal to the Dehn twist.

\hspace{20pt} Thurston's earthquake deformation is a generalization of the Fenchel--Nielsen twist. It is a deformation cutting a hyperbolic surface along a lamination instead of a simple closed geodesic (\cite{Thu}). An earthquake  deformation is one of the main actors of the two-dimensional hyperbolic geometry and used to solve the Nielsen realization problem by Kerckhoff (\cite{Ker}). Moreover, the density of the earthquake flow in the moduli space is studied by Mirzakhani (\cite{Mir}).

\hspace{20pt} In the case of closed surfaces, for any earthquake deformation $E$, there exists a sequence of the Fenchel--Nielsen twists converging $E$. In the case of marked surfaces, instead of the Fenchel--Nielsen twists, for any earthquake deformation $E$, there exists a sequence of earthquake deformations along ideal arcs converging $E$. We calculated earthquake deformations along ideal arcs (\cite{Asa} and \cite{AIK}).

\begin{figure}[H]
    \begin{center}
        \includegraphics[width=6cm]{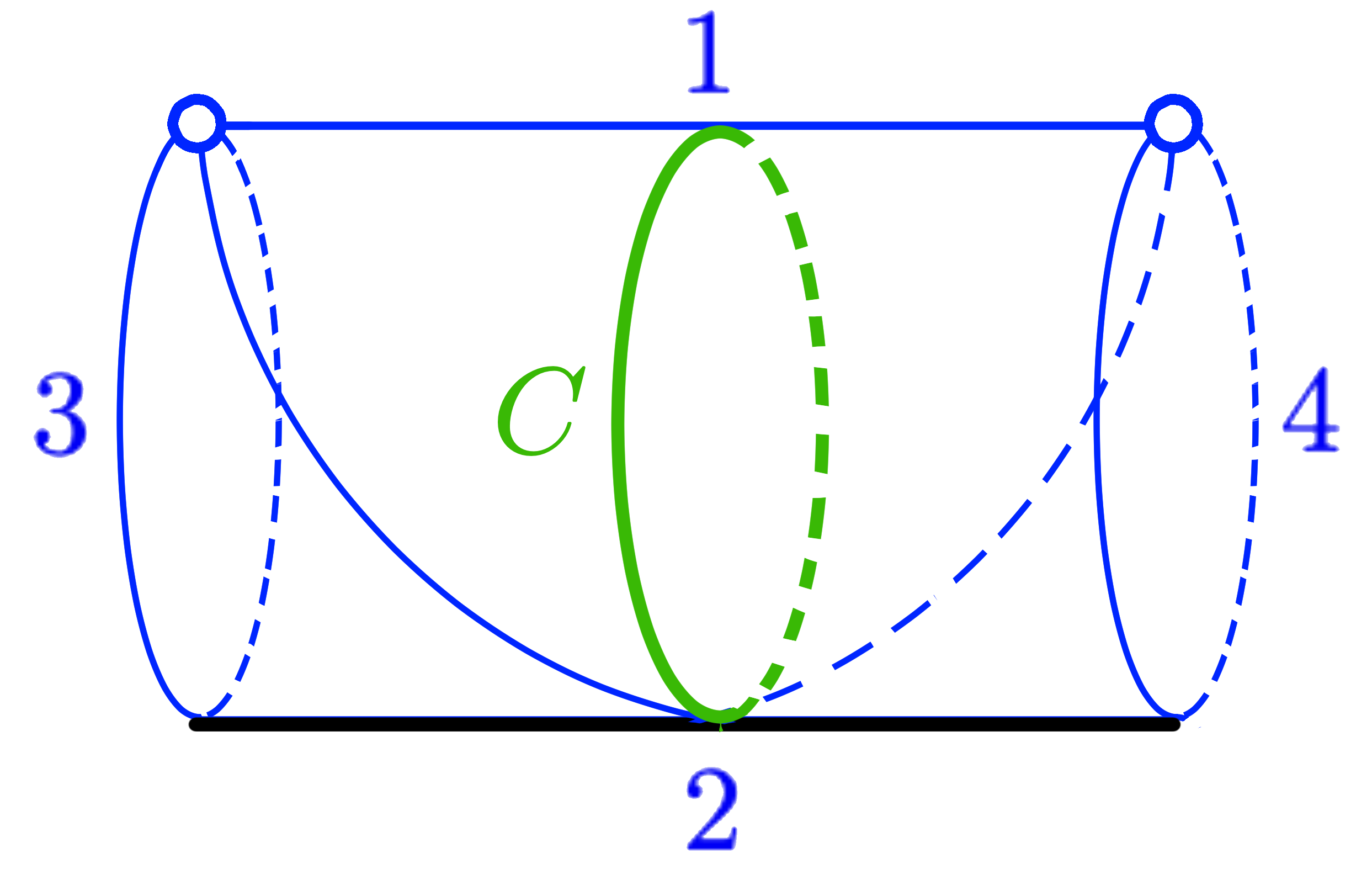}
    \end{center}
    \vspace{-30pt}
    \caption{the annulus $\annulus$}
    \label{intro}
\end{figure}

\hspace{20pt}In this paper, we calculate the Fenchel--Nielsen twist in the enhanced \Teich space of a marked surface by the cross ratio coordinates.
For a marked surface $(\Sigma, M)$ and a closed curve $C$ on $\Sigma\setminus M$, we assume that any connected component of $\Sigma\setminus C$ has at least one marked point. Then, we can embed the annulus $\annulus$ with one point on each boundary component  into $\Sigma$ so that $\annulus$ includes $C$ as Figure \ref{intro}.

\noindent We take an element $g$ of the enhanced \Teich space $\wT$ of $(\Sigma, M)$. Let the cross ratio coordinates of $g$ for the numbered ideal arcs in Figure \ref{intro} be $(X_1,X_2,X_3,X_4)$ and the length of the simple closed geodesic on $g$ corresponding to $C$ be $\Lf$. We take a rotation parameter $t\in\mathbb{R}_{\geq 0}$.
\newpage
\begin{thm}\label{main}
    In the situation above, the cross ratio coordinates $(X_1(E_{tC}(g)),X_2(E_{tC}(g)),$\\$X_3(E_{tC}(g)),$ $X_4(E_{tC}(g)))$  which the Fenchel--Nielsen twist along $C$ with length $t\Lf$ acts on are
    \begin{align*}
      X_1(E_{tC}(g))
      =&
      X_1\cdot\dfrac{2\left(X_1\,X_2\,\cosh(\Lf)-2\,\sqrt{X_1\,X_2}\,\cosh\left(\frac{\Lf}{2}\right)+X_1+1\right)\,e^{t\,\Lf}}{\left(\left(\sqrt{X_1\,X_2}\,e^{-\frac{\Lf}{2}}-X_1-1\right)\,e^{t\,\Lf}-\left(\sqrt{X_1\,X_2}\,e^{\frac{\Lf}{2}}-X_1-1\right)\right)^2}, \\
      X_2(E_{tC}(g))
      =&X_2\cdot\frac{\left(\left(\sqrt{X_1\,X_2}\,e^{-\frac{\Lf}{2}}-1\right)\,e^{t\,\Lf}-\left(\sqrt{X_1\,X_2}\,e^{\frac{\Lf}{2}}-1\right)\right)^2}{2\left(X_1\,X_2\,\cosh\,(\Lf)-2\,\sqrt{X_1\,X_2}\,\cosh\,\left(\frac{\Lf}{2}\right)+X_1+1\right)\,e^{t\,\Lf}}, \\
      X_3(E_{tC}(g))
      =&
      X_3\cdot\dfrac{\left(\sqrt{X_1\,X_2}\,e^{-\frac{\Lf}{2}}-X_1-1\right)\,e^{t\,\Lf}-\left(\sqrt{X_1\,X_2}\,e^{\frac{\Lf}{2}}-X_1-1\right)}{\left(\sqrt{X_1\,X_2}\,e^{-\frac{\Lf}{2}}-1\right)\,e^{t\,\Lf}-\left(\sqrt{X_1\,X_2}\,e^{\frac{\Lf}{2}}-1\right)}, \\
      X_4(E_{tC}(g))
      =&
      X_4\cdot\dfrac{\left(\sqrt{X_1\,X_2}\,e^{-\frac{\Lf}{2}}-X_1-1\right)\,e^{t\,\Lf}-\left(\sqrt{X_1\,X_2}\,e^{\frac{\Lf}{2}}-X_1-1\right)}{\left(\sqrt{X_1\,X_2}\,e^{-\frac{\Lf}{2}}-1\right)\,e^{t\,\Lf}-\left(\sqrt{X_1\,X_2}\,e^{\frac{\Lf}{2}}-1\right)}.
    \end{align*}
\end{thm}
\begin{rem}
    When $t$ is an integer $m$, the Fenchel--Nielsen twist becomes the $m$ times Dehn twist and $X_i(E_{tC}(g))\ (i=1,...,4)$ becomes the rational function of $X_i\  (i=1,...,4)$.
\end{rem}
\hspace{20pt}Also, the Fenchel--Nielsen twist was calculated by the coordinates of triangle length or by the trace coordinates by Garden  (\cite{Gar}).

\subsection{Motivation from the cluster algebra and future works}
\hspace{20pt}From \cite{FG2}, the cross ratio coordinates are deeply related to the cluster algebra introduced by Fomin--Zelevinsky \cite{FZ1}.
A mutation class $\bs$ is one of the subjects in the cluster algebra and associated with the piecewise-linear manifold $\tropmfd$ homeomorphic to $(\mathbb{R}^\trop)^n$. Let the fan having the cones, one of whose coordinates of $\tropmfd$ are all equal to or greater than zero as the full dimensional cones  be called the \emph{Fock--Goncharov fan} $\mathfrak{F}_\bs^+$. We denote the support of the Fock-Goncharov fan as $|\mathfrak{F}_\bs^+|$.

\begin{enumerate}
    \item The condition $|\mathfrak{F}_\bs^+| = \tropmfd$ is equivalent to the condition $\bs$ is of finite type ([FG2]). In this case, we defined a \emph{cluster earthquake map} and proved the earthquake theorem (\cite{AIK}).
    \item As for the condition $|\mathfrak{F}_\bs^+|$ is dense in $\tropmfd$, some sufficient conditions were given ([Yur1] and [Yur2]). In this case, it may be possible that we define a cluster earthquake map and prove the earthquake theorem by continuous extension of the definition on $|\mathfrak{F}_\bs^+|$.
\end{enumerate}
  
\begin{figure}[H]
    \begin{center}
        \includegraphics[width=9cm]{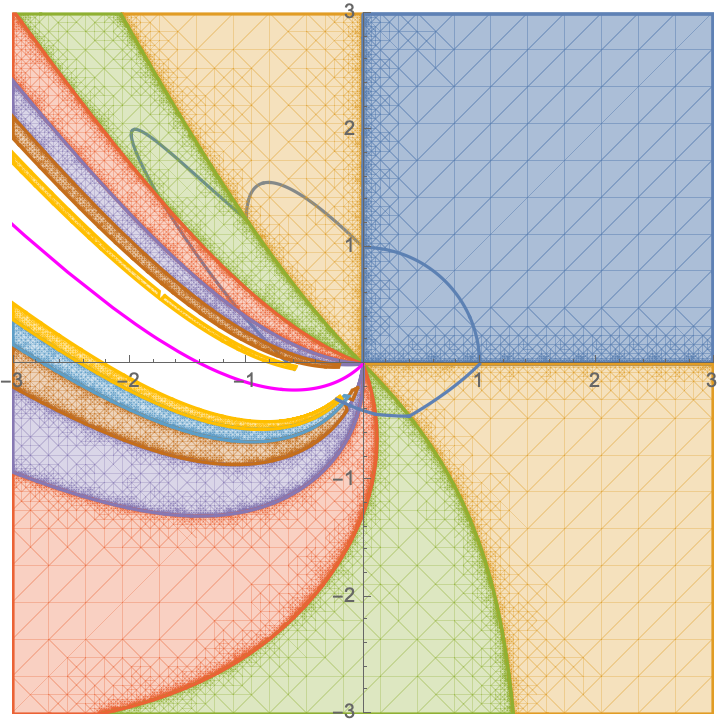}
    \end{center}
    \caption{Each color-coded area is the place where one point of $\wTone$ moves by an earthquake map for each full dimensional cone of $\mathfrak{F}_\bs^+$. The royalblue curve is the locus to which the unit circumference about some cross ratio coordinates on $\tropmfd$ is moved by an earthquake map. The magenta curve is drawn by the formula of the Fenchel--Nielsen twist in this paper. This is an incomplete figure, as it requires infinite iterations of drawing.}
    \label{zu}
\end{figure}

\hspace{20pt} The Fenchel--Nielsen twist is an earthquake on the other area than $|\mathfrak{F}_\bs^+|$. Future prospects are the explicit definition of the cluster earthquake map  and  another proof of the earthquake theorem by using the formula of the Fenchel--Nielsen twist. For example, we can draw an earthquake map on the enhanced \Teich space $\wTone$ for an annulus $\annulus$ with one point on each boundary component in Figure \ref{zu}. The marked surface $\annulus$ is the simplest case that is not treated in either \cite{BKS} or \cite{AIK}.  

\hspace{20pt}Each full dimensional cone of $\mathfrak{F}_\bs^+$ corresponds to the set of measured laminations along an ideal triangulation. We observe that the images of the earthquake map relative to $\mathfrak{F}_\bs^+$ asymptotically approach the locus of the Fenchel--Nielsen twist.

\section{The enhanced \Teich space and the cross ratio coordinates}
\hspace{20pt}We start with the same settings as \cite{AIK}. Let $\Sigma$ be an oriented
compact topological surface which may have
boundaries and $M$ be a finite subset in $\Sigma$. We assume that each boundary component has at least one point of $M$. We call $(\Sigma,M)$ a \emph{marked surface}.
  
\hspace{20pt}Let $M_\circ$ be the intersection of $M$ and the interior of $\Sigma$. We define $\Sigma^\ast:=\Sigma \setminus M_\circ$. An \emph{ideal arc} in $(\Sigma,M)$ is an isotopy class of an arc in $\Sigma^\ast$ connecting the points of $M$ without self-intersection except for its endpoints. We require that an ideal arc not to be represented by a boundary segment. A nonempty maximal disjoint collection $\Delta$ of ideal arcs without a self-folded triangle is called an \emph{ideal triangulation} of $\Sigma$. We assume that there is an ideal triangulation in $(\Sigma,M)$. A \emph{labeled triangulation} is an ideal triangulation $\Delta$ equipped with a bijection
$\ell: I \to \Delta$, where $n$ is the number of ideal arcs and
$I:=\{1,\dots,n\}$
is a fixed index set.

\begin{defi}
    Let $\mathrm{Hyp}(\Sigma)$ be the set of hyperbolic metrics on $\Sigma \setminus M$ where each point in $M\cap \Intr \Sigma$ corresponds to a funnel or a cusp, and each point in $M\cap \partial \Sigma$ corresponds to a spike of a crown. For $h \in \mathrm{Hyp}(\Sigma)$, let $M_\circ^h \subset M_\circ$ denotes the set of punctures that give rise to funnels. 
\end{defi}

\begin{figure}[H]
  \begin{center}
    \includegraphics[width=15cm]{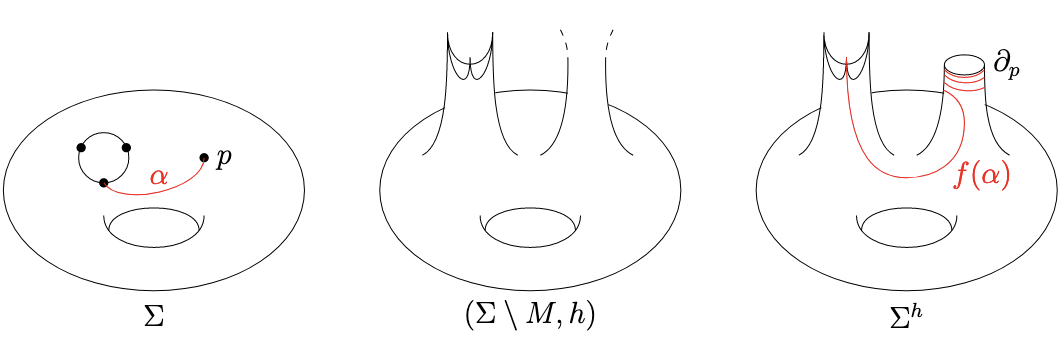}
  \end{center}
  \caption{an example of a signature}
  \label{teich}
\end{figure}

\begin{defi}
    \begin{enumerate}\label{def_Teich}
        \item For $h\in{\rm Hyp}(\Sigma)$, let $\Sigma^h$ be the hyperbolic surface obtained from $(\Sigma \setminus M,h)$ by truncating the outer side of the shortest closed geodesic in each funnel. For $p \in M_\circ^h$, let $\partial_p$ denote the resulting boundary component.
        \item  We call an orientation-preserving homeomorphism
        $f: \Intr \Sigma\rightarrow \Intr \Sigma^h$ a \emph{signed homeomorphism} if it maps a representative of each ideal arc to a complete geodesic so that each end
        incident to $p \in M_\circ^h$ is mapped to a spiraling geodesic into the geodesic boundary $\partial_p$ in either left or right direction.
        Given a signed homeomorphism $f$, its \emph{signature}
        $\epsilon_f=(\eta_p)_p\in\{+,0,-\}^{M_\circ}$ is defined by
        \[
            \eta_p=\begin{cases}
            + &\mbox{if $p\in M_\circ^h$ and $f(\alpha_p)$ spirals to the left along the geodesic boundary $\partial_p$},\\
                0&\mbox{if $p \in M_\circ \setminus M_\circ^h$}, \\
                - &\mbox{if $p\in M_\circ^h$ and $f(\alpha_p)$ spirals to the right along the geodesic boundary $\partial_p$}.
            \end{cases}
        \]
        Here $\alpha_p$ is any ideal arc incident to $p$. See also Figure \ref{teich}.
        \item We say that two pairs $(\Sigma^{h_i},f_i)$ for $i=1,2$ are \emph{equivalent} if $f_2\circ f_1^{-1}$ is homotopic to an isometry $\Sigma^{h_1} \to \Sigma^{h_2}$, and $\epsilon_{f_1}=\epsilon_{f_2}$. We call the set $\widehat{\cT}(\Sigma)$ of the equivalent classes of the pairs $(\Sigma^h,f)$ the \emph{enhanced \Teich space}.
    \end{enumerate}
\end{defi}

\hspace{20pt}We introduce the \emph{cross ratio coordinates} into  $\wT$. Let $\HHH$ be the upper half space in the complex plane and $\SSS=\mathbb{R}\cup\{\infty\}$ be the boundary of $\HHH$.

\begin{defi}
    Given a labeled triangulation $(\Delta,\ell)$ and $i \in I$, let $\square_i$ be the unique quadrilateral in $\Delta$ containing the ideal arc $\ell(i)$ as a diagonal. For a point $[\Sigma^h,f] \in \widehat{\cT}(\Sigma)$, choose a geodesic lift $\widetilde{f(\square_i)} \subset \widetilde{\Sigma^h} \subset \HHH$ of $f(\square_i)$ with respect to the universal covering $\widetilde{\Sigma^h} \to \Sigma^h$ determined by the hyperbolic structure $h$. Let $x,y,z,w \in \SSS$ be the four ideal vertices of $\widetilde{f(\square_i)}$ in this counter-clockwise order, and assume that $x$ is an endpoint of the lift of $\ell(i)$.
    Then we define the \emph{cross ratio coordinate}  associated with $i \in I$ by
    \[X_i^{(\Delta,\ell)}([\Sigma^h,f])
    \ceq [x : y : z : w] >0,
    \]
    where
    \begin{align*}
        [x:y:z:w]&\ceq \frac{w-x}{w-z}\cdot\frac{z-y}{y-x}
    \end{align*}
    denotes the cross ratio of the four points in $\SSS$.
    Since the right-hand side is invariant under M\"{o}bius transformations, $X_i^{(\Delta,\ell)}$ does not depend on the choice of the lift $\widetilde{f(\square_i)}$, and does not change if we choose the other endpoint of $\ell(i)$ as $x$.
\end{defi}
 In order to define the Fenchel--Nielsen twist, we give another description of $\wT$.
 We take $[\Sigma^h_0,f_0]\in\wT$
and denote the closure in $\HHH\cup\SSS$ of the lift $\widetilde{\Sigma^h_0}$ of $\Sigma^h_0$ as $\overline{\widetilde{\Sigma^h_0}}$. Let $\Fi=\overline{\widetilde{\Sigma^h_0}}\cap\SSS$. We consider that $\Fi$ has the same cyclic order as $\SSS$.

We consider the set
\[
  \left\{
    (\rho, \psi) \middle|
    \begin{array}{lll}
      \rho:\pi_1(\Sigma)\rightarrow\PSL \text { is a Fuchsian representation,}  \\
      \psi:\Fi\rightarrow\SSS \text{ is }\rho(\pi_1(\Sigma))\text{-equivalent and cyclic-order-preserving.}
    \end{array}
  \right\}
\]
and say that $(\rho_i,\psi_i)\ (i=1,2)$ are \emph{equivalent} if there exists $\gamma\in\PSL$ satisfying $\rho_2=\gamma\cdot\rho_1\cdot\gamma^{-1}$ and $\psi_2=\gamma\cdot\psi_1$.
We denote the set of all of these equivalent classes as $\wiT$.

\begin{prop}[cf.~{\cite[Theorem 1.6]{FG1}}]\label{Farey}
    There is a natural bijection between the enhanced \Teich\ space $\wT$ and $\wiT$.
\end{prop}
Indeed, given $[\Sigma^h,f] \in \widehat{\cT}(\Sigma)$, the representation $\rho$ is the monodromy of the hyperbolic structure $h$, and $\psi$ is given by the continuous extension of the lift $\widetilde{f}:\widetilde{\Sigma} \to \widetilde{\Sigma^h} \subset \HHH$ to $\SSS$.

\hspace{20pt}For a simple closed curve $C$ of $\Sigma\setminus M$ and $t\in\mathbb{R}_{\geq 0}$, we define \emph{the Fenchel--Nielsen twist} $E_{tC}:\wiT\rightarrow\wiT$ as follows. We take $[\rho,\psi]\in\wiT$ and let $\widetilde{\Sigma_\rho}$ be the universal cover of the hyperbolic surface $\Sigma_\rho:=\HHH/\rho(\pi_1(\Sigma))$. Let $\widetilde{C}$ be the union of all lifts of $C$ in $\widetilde{\Sigma_\rho}$. We call each connected component of $\widetilde{\Sigma_\rho}\setminus \widetilde{C}$ a \emph{gap}. We fix a gap $S_0$. Let the length of $C$ in $\Sigma_\rho$ be $L(C)$. We endow each gap $S$ with a hyperbolic element $E_S\in\PSL$ as follows.
\begin{itemize}
  \item $E_{S_0}$ is the identity.
  \item For any gap $S_1$ adjacent to $S_0$, let $l_1$ be the intervening geodesic between these gaps. Let $E_{S_1}$ be the hyperbolic element whose axis is $l_1$, whose translation distance is $t\cdot L(C)$ and whose direction is to the left (i.e. the restriction $E_{S_1}|_{S_1}$ of $E_{S_1}$ moves $S_1$ to the left from $S_0$'s perspective). 
  \item For any gap $S_2$ adjacent to $S_1$ and other than $S_0$, let $l_2$ be the intervening geodesic between these gaps. Let $\overline{E_{S_2}}$ be the hyperbolic element whose axis is $l_2$, whose translation distance is $t\cdot L(C)$ and whose direction is to the left. We define $E_{S_2}:=E_{S_1}\circ\overline{E_{S_2}}$.
  \item Inductively iterating, for any gap $S_n$ adjacent to $S_{n-1}$ and other than $S_{n-2}$, let $l_n$ be the intervening geodesic between these gaps. Let $\overline{E_{S_n}}$ be the hyperbolic element whose axis is $l_n$, whose translation distance is $t\cdot L(C)$ and whose direction is to the left. We define $E_{S_n}:=E_{S_{n-1}}\circ \overline{E_{S_n}}$.
\end{itemize}
We define $\widetilde{E_{tC}}:\widetilde{\Sigma_\rho}\setminus \widetilde{C}\to \HHH$ as $\widetilde{E_{tC}}:=\displaystyle\bigcup_{S\text{ : gap}}E_S|_S$ and continuously extend it to $\partial_\infty\widetilde{E_{tC}}:\partial_\infty\widetilde{\Sigma_\rho}\rightarrow\SSS$, where $\partial_\infty\widetilde{\Sigma_\rho}$ is the intersection of $\widetilde{\Sigma_\rho}$ and $\SSS$.
\begin{enumerate}
   \item We define a group-homomorphism $\rho':\pi_1(\Sigma)\to \PSL$ by the condition $\partial_\infty \widetilde{E_{tC}}\circ\rho(\gamma)=\rho'(\gamma)\circ\partial_\infty \widetilde{E_{tC}}$ on $\partial_\infty\widetilde{\Sigma_\rho}$ for any $\gamma\in\pi_1(\Sigma)$. Then, $\rho'$ is a Fuchsian representation by \cite[Proposition 3.2.]{BKS}.
   \item Define $\psi':\cF_\infty(\Sigma) \to \SSS$ by $\psi':=\partial_\infty \widetilde{E_{tC}}\circ \psi$. Then, $\psi'$ is $\rho'$-equivariant and order-preserving. See \cite{BB,BKS}.
\end{enumerate}
As above, we define
\[
 E_{tC}:\wiT\rightarrow\wiT,\ [\rho,\psi]\mapsto [\rho',\psi'].
\]
Moreover, $E_{tC}:\wT\rightarrow\wT$ is naturally induced by Proposition \ref{Farey}.

\begin{defi}
    We call $E_{tC}:\wT\rightarrow\wT$ the \emph{Fenchel--Nielsen twist} along a simple closed curve $C$ with a twist parameter $t$.
\end{defi}

\section{A formula of the Fenchel--Nielsen twist}
\subsection{The proof of theorem \ref{main}}
\hspace{20pt}First, we consider an annulus $\annulus$ with one point on each boundary component and ideal arcs 
in Figure \ref{intro}. We take $g\in\wTone$ and calculate the Fenchel--Nielsen twist along the 
non-trivial and non-self-intersecting closed curve $C$. We number each edge and let the cross ratio coordinates of $g$ be $(X_1, X_2, X_3, X_4)$.

\begin{figure}[H]
    \begin{center}
        \includegraphics[width=7cm]{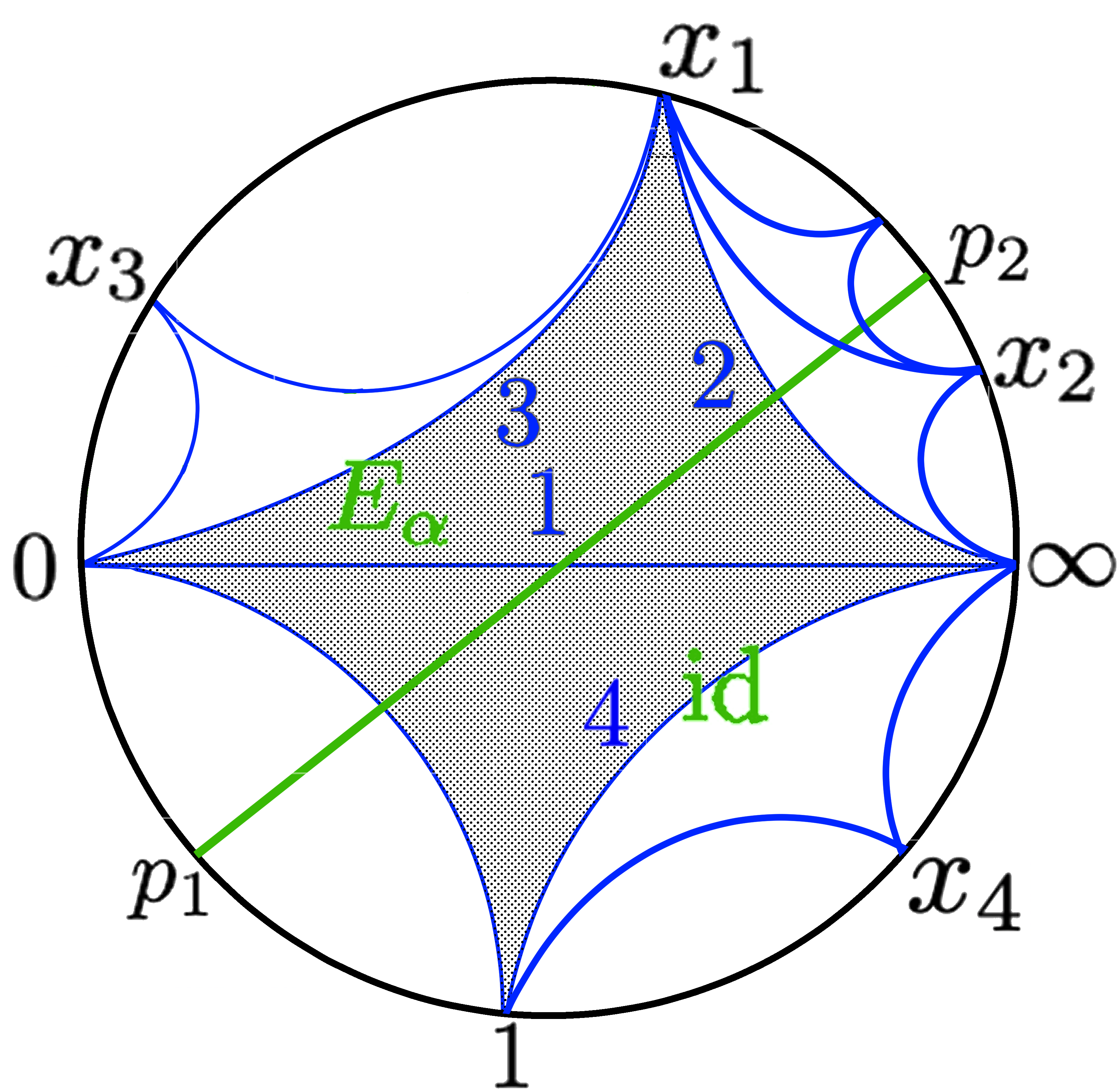}
    \end{center}
    \vspace{-25pt}
    \caption{the universal cover $\HHH$}
    \label{lift}
\end{figure}

We lift the surfaces to $\HHH$ and draw the fundamental domain of the surfaces with dots. Let the endpoints of the lift of the ideal arcs in $\SSS$ be $x_1,\ x_2,\ 0,\ x_3,\ 1,\ x_4\ \infty$ in Figure \ref{lift}.

Then, we get
\begin{align*}
  \ x_1 = -X_{1},\ x_2 = -X_{1}\,\left(X_{2}+1\right),\ x_3 = -\dfrac{X_{1}\,X_{3}}{X_{3}+1},\ x_4 = \dfrac{X_{4}+1}{X_{4}}.
\end{align*}

We denote the hyperbolic element attaching the edge $2$ as $f_2$. The element $f_2$ maps $0$, $1$ and $\infty$ to $x_1,\ \infty$ and $x_2$, respectively. Here,
\[
  f_2 = \dfrac{1}{\sqrt{X_1X_2}} \begin{bmatrix}
	  X_{1}\,(X_{2}+1) & -X_{1}\\
     -1 & 1
 \end{bmatrix}\in \PSL
\]
and we calculate the absolute value of the trace of $f_2$
\[
  \trf = \dfrac{X_1(X_2 + 1) + 1}{\sqrt{X_1X_2}}
\]
and the hyperbolic distance of $f_2$
\[
  L(f_2) = 2\,\cosh^{-1}\left(\frac{\trf}{2}\right).
\]

The fixed points of $f_2$ are
\begin{align*}
    p_1&:= \frac{-X_{1}\,(X_{2}+1)+1+\sqrt{{(X_{1}\,(X_{2}+1)+1)}^2-4\,X_{1}\,X_{2}}}{2}>0, \\
    p_2 &:=\frac{-X_{1}\,(X_{2}+1)+1-\sqrt{{(X_{1}\,(X_{2}+1)+1)}^2-4\,X_{1}\,X_{2}}}{2}<0.
\end{align*}

We normalize the scale of the Fenchel--Nielsen deformation by $L(f_2)$. Therefore, We let $t\ (t \geq 0)$ be a twist parameter and define a stratum map $E_\alpha$ as
\begin{align*}
A_1&:= \dfrac{1}{\sqrt{p_1-p_2}} \begin{bmatrix}
1 & -p_1 \\
1 & -p_2 \\
\end{bmatrix},\\
E_{\alpha}&:= A_1^{-1}\cdot \begin{bmatrix}
e^{\frac{t\,L(f_2)}{2}} & 0 \\
0 & e^{-\frac{t\,L(f_2)}{2}}
\end{bmatrix} \cdot A_1\\
&=\dfrac{1}{\left(p_1-p_2\right)e^{\frac{t\,L(f_2)}{2}}}\left[\begin{array}{cc}
p_1-p_2\,e^{t\,\Lf} & p_1\,p_2\,(e^{t\,\Lf}-1) \\
-e^{t\,\Lf}+1 & p_1e^{t\,\Lf}-p_2
\end{array}\right].
\end{align*}

The stratum maps related to cross ratio coordinates are the identity and $E_{\alpha}$ as Figure \ref{lift}.

The values acted by the lift $\widetilde{E_{tC}}$ of $E_{tC}$ are
\begin{align*}
\widetilde{E_{tC}}(0)&=\dfrac{p_1\,p_2(t-1)}{p_1\,e^{t\,\Lf}-p_2},\\
\widetilde{E_{tC}}(x_1)&=\dfrac{(X_1\,p_2+p_1\,p_2)\,e^{t\,\Lf}-(X_1\,p_1+p_1\,p_2)}{(X_1+p_1)\,e^{t\,\Lf}-(X_1+p_2)},\\
\widetilde{E_{tC}}(x_3)&=\dfrac{(X_1\,X_3\,p_2+(X_3+1)\,p_1\,p_2)\,e^{t\,\Lf}-(X_1\,X_3\,p_1+(X_3+1)\,p_1\,p_2)}{(X_1\,X_3+(X_3+1)\,p_1)\,e^{t\,\Lf}-(X_1\,X_3+(X_3+1)\,p_2)}.
\end{align*}

When we use
\[
  p_1\,p_2 = -X_1,
\]
we get
\begin{align*}
\widetilde{E_{tC}}(0)&=-\dfrac{X_1\,(t-1)}{p_1\,e^{t\,\Lf}-p_2},\\
\widetilde{E_{tC}}(x_1)&=\dfrac{X_1\,((p_2-1)\,e^{t\,\Lf}-(p_1-1))}{(X_1+p_1)\,e^{t\,\Lf}-(X_1+p_2)},\\
\widetilde{E_{tC}}(x_3)&=\dfrac{X_1\,((X_3\,p_2-X_3-1)\,e^{t\,\Lf}-(X_3\,p_1-X_3-1))}{(X_1\,X_3+(X_3+1)\,p_1)\,e^{t\,\Lf}-(X_1\,X_3+(X_3+1)\,p_2)}.
\end{align*}

We calculate the cross ratio coordinates $(X_1(E_{tC}(g)), X_2(E_{tC}(g)), X_3(E_{tC}(g)),X_4(E_{tC}(g)))$ which  $E_{tC}$ acts on.

\begin{align*}
X_1(E_{tC}(g))=
\hspace{15pt}[E_\alpha\left(0\right):1:\infty:E_\alpha\left(X_1\right)]\hspace{15pt}=
X_1\cdot\dfrac{(p_1^2+p_2^2+2\,X_1)\,e^{t\,\Lf}}{((X_1+p_1)\,e^{t\,\Lf}-(X_1+p_2))^2},\hspace{58pt}\\
X_2(E_{tC}(g))=
[E_\alpha\left(0\right):E_\alpha\left(X_2\right):1:\infty]\hspace{290pt}\\
=\dfrac{\left(p_1 t-p_2\right)\left(\left(\left(X_2+1\right) p_1+p_2+X_1 X_2+X_1-1\right) \,e^{t\,\Lf}-\left(p_1+\left(X_2+1\right) p_2+X_1 X_2+X_1-1\right)\right)}{\left(p_1^2+p_2^2+2 X_1\right)\,e^{t\,\Lf}},\\
X_3(E_{tC}(g))=
[E_\alpha\left(0\right):\infty:E_\alpha\left(X_1\right):E_\alpha\left(X_3\right)]=
X_3\cdot\dfrac{(X_1+p_1)\,e^{t\,\Lf}-(X_1+p_2)}{p_1\,e^{t\,\Lf}-p_2},\hspace{70pt}\\
X_4(E_{tC}(g))=
\hspace{30pt}[1:X_4:\infty:E_\alpha\left(0\right)]\hspace{30pt}=
X_4\cdot\dfrac{(X_1+p_1)\,e^{t\,\Lf}-(X_1+p_2)}{p_1\,e^{t\,\Lf}-p_2}.\hspace{66pt}
\end{align*}

When we use
\begin{align*}
  &p_1 = \dfrac{1}{2}(-\sqrt{X_1\,X_2}\,\trf + 2 + \sqrt{X_1\,X_2\,(\trf-4)}))=1-\sqrt{X_1\,X_2}\,\exp \left(-\frac{\Lf}{2}\right),\\
  &p_2 = \dfrac{1}{2}(-\sqrt{X_1\,X_2}\,\trf + 2 - \sqrt{X_1\,X_2\,(\trf-4)}))=1-\sqrt{X_1\,X_2}\,\exp \left(\frac{\Lf}{2}\right),\\
  &p_1+p_2 = -X_1\,X_2 - X_1 + 1,
\end{align*}

we get
\begin{align*}
    X_1(E_{tC}(g))
    =&
    X_1\cdot\dfrac{2\left(X_1\,X_2\,\cosh(\Lf)-2\,\sqrt{X_1\,X_2}\,\cosh\left(\frac{\Lf}{2}\right)+X_1+1\right)\,e^{t\,\Lf}}{\left(\left(\sqrt{X_1\,X_2}\,e^{-\frac{\Lf}{2}}-X_1-1\right)\,e^{t\,\Lf}-\left(\sqrt{X_1\,X_2}\,e^{\frac{\Lf}{2}}-X_1-1\right)\right)^2}, \\
    X_2(E_{tC}(g))
    =&X_2\cdot\frac{\left(\left(\sqrt{X_1\,X_2}\,e^{-\frac{\Lf}{2}}-1\right)\,e^{t\,\Lf}-\left(\sqrt{X_1\,X_2}\,e^{\frac{\Lf}{2}}-1\right)\right)^2}{2\left(X_1\,X_2\,\cosh\,(\Lf)-2\,\sqrt{X_1\,X_2}\,\cosh\,\left(\frac{\Lf}{2}\right)+X_1+1\right)\,e^{t\,\Lf}},\\
    X_3(E_{tC}(g))
    =&
    X_3\cdot\dfrac{\left(\sqrt{X_1\,X_2}\,e^{-\frac{\Lf}{2}}-X_1-1\right)\,e^{t\,\Lf}-\left(\sqrt{X_1\,X_2}\,e^{\frac{\Lf}{2}}-X_1-1\right)}{\left(\sqrt{X_1\,X_2}\,e^{-\frac{\Lf}{2}}-1\right)\,e^{t\,\Lf}-\left(\sqrt{X_1\,X_2}\,e^{\frac{\Lf}{2}}-1\right)}, \\
    X_4(E_{tC}(g))
    =&
    X_4\cdot\dfrac{\left(\sqrt{X_1\,X_2}\,e^{-\frac{\Lf}{2}}-X_1-1\right)\,e^{t\,\Lf}-\left(\sqrt{X_1\,X_2}\,e^{\frac{\Lf}{2}}-X_1-1\right)}{\left(\sqrt{X_1\,X_2}\,e^{-\frac{\Lf}{2}}-1\right)\,e^{t\,\Lf}-\left(\sqrt{X_1\,X_2}\,e^{\frac{\Lf}{2}}-1\right)}.
    \end{align*}
    
    We have gotten the formula for $\annulus$.

    \hspace{20pt}Finally, we consider a general case. We take a marked surface $(\Sigma,M)$ and a simple closed curve $C$ on $\Sigma\setminus M$. When any connected component of $\Sigma\setminus C$ has at least one marked point, there are two ideal arcs constituting the boundary of $\annulus$ which includes $C$ and whose internal has no marked point. Then, we introduce an ideal triangulation into $\annulus$ in the same way as Figure \ref{intro}, add ideal arcs and produce an ideal triangulation of $(\Sigma,M)$. We have proved Theorem \ref{main}.

    \subsection{The Dehn twist}
    We consider the formula of the Fenchel--Nielsen twist in the case of $t=0$ and $t=1$. By the trace and hyperbolic distance of $f_2$, we get
    \begin{align*}
      &\cosh \left(
        \dfrac{\Lf}{2}
      \right)
      =
      \dfrac{1}{2}\left(
        \dfrac{X_1\,X_2 + X_1 + 1}{\sqrt{X_1\,X_2}}
      \right),\\
      &\cosh (\Lf)
      =
      \dfrac{(X_1\,X_2 + X_1 + 1)^2 - 2\,X_1\,X_2}{2\,X_1\,X_2}.
    \end{align*}
    \begin{enumerate}
      \item When $t=0$,
        \begin{align*}
        X_1(E_{0\cdot C}(g))
        &=
        X_1\,\cdot\,\frac{2\left(X_1\,X_2\,\cosh\,\left(L\left(f_2\right)\right)-2\,\sqrt{X_1\,X_2}\,\cosh\,\left(\frac{L\left(f_2\right)}{2}\right)+X_1+1\right)}{2\,X_1\,X_2\left(\cosh\,\left(L\left(f_2\right)\right)-1\right)}\\
        &=X_1.
        \end{align*}
        The other coordinates are calculated in the same way and we get
        \[
          (X_1(E_{0\cdot C}(g)),X_2(E_{0\cdot C}(g)),X_3(E_{0\cdot C}(g)),X_4(E_{0\cdot C}(g)))=(X_1,X_2,X_3,X_4).
        \]
      \item When $t=1$,
      \begin{align*}
      X_1(E_{1\cdot C}(g))
      &=
      X_1\,\cdot\,\frac{X_1\,X_2\,\cosh\,\left(L\left(f_2\right)\right)-2\,\sqrt{X_1\,X_2}\,\cosh\,\left(\frac{L\left(f_2\right)}{2}\right)+X_1+1}{\left(X_1+1\right)^2\left(\cosh\,\left(L\left(f_2\right)\right)-1\right)}\\
      =\dfrac{X_1^2\,X_2}{(X_1 + 1)^2}.
      \end{align*}
      The other coordinates are calculated in the same way and we get
      \begin{align*}
        &(X_1(E_{1\cdot C}(g)),X_2(E_{1\cdot C}(g)),X_3(E_{1\cdot C}(g)),X_4(E_{1\cdot C}(g)))\\
        &=\left(
          \dfrac{X_1^2\, X_2}{(X_1 + 1)^2},
          \dfrac{1}{X_1},
          (X_1 + 1)\,X_3,
          (X_1 + 1)\,X_4
        \right).
      \end{align*}
      It is a formula of the Dehn twist.
    \end{enumerate}

    \section*{Acknowledgement}
    The author thanks Tsukasa Ishibashi and Shunsuke Kano for their valuable discussion. He would like to thank his supervisor, Takuya Sakasai for his patient guidance.

    \bibliography{infinite_type}
    \bibliographystyle{plain}
    

    \vspace{30pt}
    Takeru Asaka
    
    Graduate School of Information Sciences, Tohoku University, 6-3 Aoba, Aramaki, Aoba-ku, Sendai, Miyagi 980-8578, Japan.
  
    Email : asakatakeru@gmail.com

\end{document}